\documentclass{amsart}

\usepackage{graphicx}
\usepackage{color}
\usepackage{eepic}
\usepackage{amssymb}
\usepackage{amsthm}
\usepackage{amsmath}
\usepackage{cite}
\usepackage{psfrag}

\newtheorem{theorem}{Theorem}%[section]

\newtheorem{conjecture}{Conjecture}

\theoremstyle{definition}
\newtheorem{example}{Example}

\newcommand{\Z}{\mathbb{Z}}
\newcommand{\Q}{\mathbb{Q}}
\newcommand{\F}{\mathbb{F}}

             \def \N {{\mathbb N}}

\def\F{\mathbb F}

\def\R{\mathbb R}
\def\C{\mathbb C}
\def\Q{\mathbb Q}
\def\Z{\mathbb Z}

\def\det{\mbox{det}}

\renewcommand{\P}{{\mathbf{P}}}
\newcommand{\GL}{{\rm GL}}
\newcommand{\GSp}{{\rm GSp}}
\newcommand{\Sp}{{\rm Sp}}
\newcommand{\im}{{\rm im\ }}
\newcommand{\mmker}{{\rm ker\ }}
\newcommand{\markmin}{{\rm min}}
\newcommand{\rank}{{\rm rank\ }}
\newcommand{\SL}{{\rm SL}}
\newcommand{\SO}{{\rm SO}}
\newcommand{\shh}{{\sc Sheafhom}}
\newcommand{\cochains}{C} % maybe $E_1^{i,j}$ later

\def\BorelSerre{\text{BS}}
\def\Eis{\text{Eis}}
\def\cusp{\text{cusp}}
\def\Frob{\text{Frob}}
\def\Gal{\text{Gal}}
\def\IIa{\text{IIa}}
\def\IIb{\text{IIb}}
\def\IIIa{\text{IIIa}}
\def\IIIb{\text{IIIb}}
\def\IV{\text{IV}}
\def\g{\text{G}}
\def\ng{\text{nG}}

\begin{document}

%%%from my template
\title{Cohomology of Congruence Subgroups of $\SL_4(\Z)$. III}

\author{Avner Ash} \address{Boston College\\ Chestnut Hill, MA 02445}
\email{Avner.Ash@bc.edu} \author{ Paul E. Gunnells}
\address{University of Massachusetts Amherst\\ Amherst, MA 01003}
\email{gunnells@math.umass.edu} \author{Mark McConnell}
\address{Center for Communications Research\\ Princeton, New Jersey 08540}
\email{mwmccon@idaccr.org}

\thanks{AA wishes to thank the National Science Foundation for support
of this research through NSF grant DMS-0455240. PG wishes to thank the
National Science Foundation for support of this research through NSF
grant DMS-0801214. We thank T.~Ibukiyama and C.~Poor for helpful
conversations.  Finally we thank the referees for
helpful references and comments.}



\begin{abstract}
In two previous papers \cite{AGM1, AGM2} we computed cohomology groups
$H^{5}(\Gamma_{0} (N); \C)$ for a range of levels~$N$, where
$\Gamma_{0} (N)$ is the congruence subgroup of $\SL_{4} (\Z)$
consisting of all matrices with bottom row congruent to $(0,0,0,*)$
mod~$N$.  In this note we update this earlier work by carrying it out
for prime levels up to $N = 211$. This requires new methods in sparse matrix
reduction, which are the main focus of the paper.  Our computations
involve matrices with up to 20~million nonzero entries.  We also make
two conjectures concerning the contributions to $H^{5}(\Gamma_{0} (N);
\C)$ for $N $ prime coming from Eisenstein series and Siegel modular
forms.
\end{abstract}

\maketitle
\section{Introduction}\label{intro}

In two previous papers \cite{AGM1, AGM2} we computed the cohomology in
degree~5 of congruence subgroups $\Gamma_0(N) \subset \SL_{4}(\Z)$
with trivial $\C$ coefficients, where $\Gamma_{0} (N)$ is the subgroup
of $\SL_{4} (\Z)$ consisting of all matrices with bottom row congruent
to $(0,0,0,*)$ mod~$N$.  The highest level we reached was $N=83$.  We
also computed some Hecke operators on these cohomology groups and
identified the cohomology as either cuspidal or as coming from the
boundary (Eisensteinian).

In this paper we concentrate on explaining new techniques we have
developed to reduce very large sparse matrices.  These techniques have
enabled us to carry out our computations for much higher values of the
level~$N$.  We explain in Section \ref{s:marks-section-top-level} that
our algorithms differ from others in the literature because we must
compute change of basis matrices.  As an oversimplified illustration,
imagine solving $A \mathbf{x} = \mathbf{b}$ for an invertible
matrix~$A$.  Classical dense methods produce an invertible change of
basis matrix~$P$ where $PA$ has a desirable form, and we solve
for~$\mathbf{x}$ by computing $P \mathbf{b}$.  When~$A$ is large and
sparse, computing~$P$ is much too expensive if finding~$\mathbf{x}$ is
our only goal.  Iterative methods like Wiedemann's
produce~$\mathbf{x}$ more simply.  In this paper, however, we compute
explicit cocycles in our cohomology groups, and compute their images
under Hecke operators.  As explained in~(\ref{subsec-coh}), change of
basis matrices are essential for this task.  (See Section
\ref{s:marks-section-top-level} for references.  The illustration is
oversimplified because the actual~$A$ have less than full rank and are
not even square.)

The linear algebra issues are compounded when computing the Smith
normal form (SNF) of integer matrices.  Modern SNF methods reduce
mod~$p$ for a range of~$p$, solve at each~$p$ by iterative methods,
and lift to a final result over~$\Z$ by Chinese remaindering
techniques.  It seems unknown, though, how to find the SNF change of
basis matrices by Chinese remaindering.  Hence we use a different
approach~(\ref{subsec-SNF}).  See~(\ref{speed-space}) for a comparison
of times.  Although iterative methods can be
efficiently parallelized, this paper does not use parallel techniques.

In future installments of this project we will look at torsion classes
in $H^5(\Gamma_0(N);\Z)$ as well as twisted coefficient modules.  For
the torsion classes, we will test Conjecture B of \cite{A} that
asserts that they have attached Galois representations.  The new
sparse matrix techniques discussed here will be of great importance in
carrying this project forward.

In the paper \cite{dutch} of Van Geemen, Van der Kallen, Top, and
Verberkmoes, cohomology classes for $\GL(3)$ were found by working
modulo small primes and using the LLL algorithm to reconstruct
integral solutions.  This is a useful method that various people
(including ourselves in the past) have followed.  However, in this
paper we work solely modulo a 5 digit prime $p$ without lifting to $\Z
$.  Lifting to $\Z $ would be a huge computatonal effort at larger
levels.  The prime $p$ is small enough to make computation fast, and
large enough to make us morally certain that we are actually finding
the complex betti numbers and Hecke eigenvalues.  The fact that all
our data is accounted for by Eisenstein series and liftings of
automorphic forms confirms this.

We continue to find that the cuspidal part consists of functorial
lifts of Siegel modular forms from paramodular subgroups of $\Sp_4 (\Q
)$ that are not Gritsenko lifts, as described in \cite{AGM2} for
levels $N=61, 73, 79$.  We conjecture that these functorial lifts will
always occur, at least for prime level, in Conjecture
\ref{conj:siegel} of Section \ref{s:conjectures.computational}.  These
lifted forms correspond to selfdual automorphic representations on
$\GL(4)/\Q$.  We were hoping to find non-lifted cuspidal cohomology
classes, which would correspond to non-selfdual automorphic
representations.  Unfortunately, we found none.  We see no reason why
they should not exist for larger~$N$, but no one has proven their
existence.  It should be noted that non-selfdual Galois
representations, that by Langlands' philosophy would predict
non-selfdual automorphic representations for $\GL(4)$ of the type we
are searching for, were constructed by J. Scholten \cite{sholten}.

Our data for the boundary cohomology for prime level leads to our
Conjecture \ref{conj:eisenstein} of Section \ref{s:conjectures.computational}
that identifies its constituents as various Eisensteinian lifts of certain
classical cuspforms of weights 2 and 4, and of certain cuspidal
cohomology classes for $\GL(3)/\Q$.  This conjecture is a refinement
of Conjecture 1 of \cite{AGM2}.

We ought to have similar conjectures for composite level, but we don't
have enough data to justify an attempt to make them.  The size of the
matrices and the complexity of computing the Hecke operators increases
as $N$ grows larger or more composite.  Therefore at a certain point
we stopped computing for composite~$N$ but were able to continue for
prime~$N$ up to level 211.  Similarly the size of the computation of
the Hecke operators at a prime $l$ increases dramatically with $l$, so
that in fact for the new levels in this paper, we compute the Hecke
operators at most for $l = 2$.

The index of $\Gamma_0(N)$ in $\SL_{4}(\Z)$ grows like $O(N^3)$.  Thus
the matrices we need to reduce are on the order of $N^3\times N^3$.
This growth in size is what makes this computational project so much
more difficult to carry out for large~$N$, compared to the cohomology
of congruence subgroups of $\SL_{2}(\Z)$, $\SL_{3}(\Z)$,
$\Sp_{4}(\Z)$, and other lower rank groups.  The implied constant in
the $O(N^3)$ is larger when~$N$ is composite, which is why we
eventually had to restrict ourselves to prime~$N$.  Also the Hecke
operator computations become much more lengthy for composite~$N$.

Please refer to \cite{AGM1} for explanations of why we look in
degree~5, the decomposition of the cohomology into the cuspidal part
and the part coming from the Borel--Serre boundary, and how we perform
the computations.  We also review there the number theoretic reasons
for being interested in this cohomology, primarily because of
connections to Galois representations and motives.  In \cite{AGM2} the
reader will find how we identified the cuspidal part of the cohomology
as lifted from $\GSp(4)/\Q$ and why this must be so for selfdual
classes of prime level.

In \cite{AGM1} we explained that the well-rounded retract $W$ of the
symmetric space for $\SL_{4}(\R)$ is a contractible cell complex on
which $\SL_{4}(\Z)$ acts cellularly with finite stabilizers of orders
divisible only by 2, 3 and~5, and that~$W$ has only finitely many
cells modulo $\SL_{4}(\Z)$.  Therefore we can use the chain complex
$C^*(W/\Gamma_0(N);\F)$ to compute $H^*(\Gamma_0(N); \F)$ for any
field~$\F$ of characteristic not equal to 2, 3 or~5.  In practice we
substitute a large finite field for $\C$ as justified in \cite{AGM1}.

Also in \cite{AGM1} we described explicitly and in detail how to
handle the data structures needed to construct the chain complex from
$W/\Gamma_0(N)$ and hence to create the matrices whose nullspaces
modulo column spaces are isomorphic to the cohomology.

In this paper we continue to use this set-up.  The new thing here is
the method explained in Section \ref{s:marks-section-top-level}, which
enables us to take~$N$ as far as 211.

In Section \ref{s:eisenstein.paramodular} we give the background on
Eisenstein cohomology and Siegel modular forms needed to present our
computational results and to formulate our conjectures.  Finally, in
Section \ref{s:conjectures.computational} we state two conjectures
about the structure of $H^{5} (\Gamma_{0} (N); \C)$ for~$N$ prime,
give the results of our computations of $H^5(\Gamma_0(N);\C)$ for~$N$
prime, $83\le N \le 211$, and verify the two conjectures in this
range.  The first, Conjecture \ref{conj:eisenstein}, improves on
\cite[Conjecture 1]{AGM1} by fixing the weight 4 part of the
Eisensteinian cohomology to those weight 4 cuspforms $f$ whose central
special value vanishes.  We also feel confident now of conjecturing
that our list of classes for the boundary cohomology is complete in
this case.  The second, Conjecture \ref{conj:siegel}, states exactly
which cuspidal paramodular Siegel forms at prime level show up in the
cuspidal cohomology.

\section{Computational Methods}
\label{s:marks-section-top-level}
Our problem is to find $H^5$ of a complex of free
$R$-modules for some ring $R$,
\begin{equation}
0 \longleftarrow \cochains^6
  \stackrel{d^5}\longleftarrow \cochains^5
  \stackrel{d^4}\longleftarrow \cochains^4
  \longleftarrow \cdots
  \label{main-complex}
\end{equation}
Let $n_i = \rank \cochains^i$.  View the~$\cochains^i$ as a space of
column vectors, and represent the~$d^i$ as matrices.  All the matrix
entries lie in~$\Z$.  It is possible to carry out our computations
over~$R = \Z$, obtaining the torsion in the cohomology along with the
free part.  Our next paper (IV in this series) will study the torsion.
In the present paper, we work over a field.  However, we will
intersperse remarks on the computational problem for more general~$R$,
with an eye to future papers in the series.

In principle we want to work over $R = \C$, because our purpose is to study automorphic forms.  To avoid
round-off error and numerical instability, however, we replace~$\C$
with a finite field~$\F_p$ of prime order~$p$.  If~$p$ is large
enough, it is extremely likely that $\dim H^i(\Gamma_0(N); \F_p)$ will
equal $\dim H^i(\Gamma_0(N); \C)$ for the~$N$ we consider, and that Hecke
computations will work compatibly.  We generally use $p = 12379$, the
fourth prime after 12345.  We give details about the choice of~$p$
in~(\ref{sparse-elt}).

\begin{figure}[htb]
\begin{center}
%\psfrag{N^3/96}{$N^{3}/96$}
%\psfrag{N^3/10}{$N^{3}/10$}
%\psfrag{25/72 N^3}{$25/72 N^{3}$}
%\includegraphics[width=3.2916667in]{matsizes2.png} % 316px / 96dpi
%\includegraphics[scale=0.7]{matrices}
%\input{matrices.pdf_t}
\input{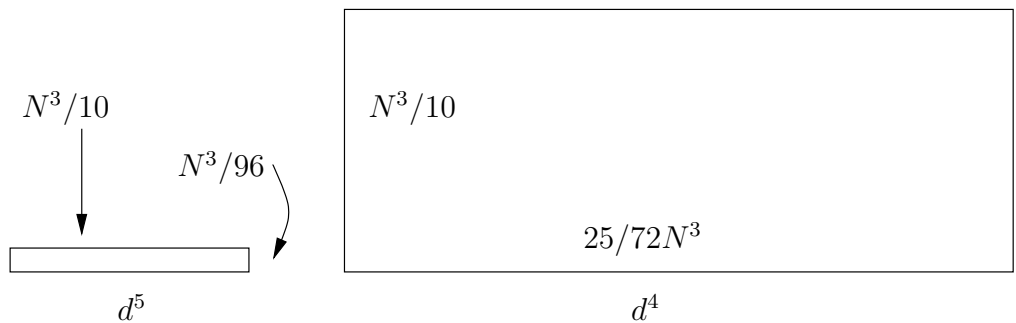}
\caption{Relative sizes of the matrices $d^5$, $d^4$ in our cochain
complex.  We will exploit the fact that $d^4$ is wider than its
height.}
\label{fig-d4-d5}
\end{center}
\end{figure}

The matrices are sparse matrices, meaning only a small fraction of the
entries in each row and column are nonzero.  Our largest matrices,
those for~$d^4$, have dimension $n_5 \times n_4 \approx N^3/10 \times
25N^3/72$ for prime~$N$.  However, at most~6 of the entries in each
column are nonzero, and at most~26 in each row.  The~6 and~26 are
independent of~$N$.  The matrices for~$d^5$ have dimension $n_6 \times
n_5 \approx N^3/96 \times N^3/10$ for prime~$N$.  All these estimates
are a few times larger for composite~$N$.  We give more precise
estimates for the~$n_i$ in~(\ref{n-i-estimates}).  The relative sizes
are shown in Figure~\ref{fig-d4-d5}.

Given an $m\times n$ sparse matrix~$A$, our fundamental computational
task is to find the Smith normal form (SNF) decomposition
\begin{equation}
A = PDQ. \label{def-SNF}
\end{equation}
Here $P\in \GL_m(R)$, $Q\in \GL_n(R)$, $D$ is zero except on the
diagonal, and the nonzero diagonal entries satisfy $d_{ii} \mid
d_{i+1,i+1}$ for all~$i$.  There is a~$\rho$ with $d_{ii} \ne 0$ for
$0 \leq i < \rho$ and $d_{ii} = 0$ for $i \geq \rho$; when~$R$ is a field,
$\rho$ is the rank of~$A$.  We call~$P$ and~$Q$ {\em change of
basis\/} matrices.

To carry out the calculations, we used \shh~2.2, a free software
package written by the third author~\cite{shh,shh-lispwire}.  \shh\
performs large-scale computations in the category of
finitely-generated $R$-modules, where~$R$ is any principal ideal
domain supporting exact computation.  Most development work in \shh\
has been for domains that are not fields, especially~$\Z$ and other
rings of integers of class number~1.  In this sense it differs from
most of the sparse matrix literature, which looks at~$\R$ and~$\C$
\cite{DER,GVL,Davis,Matlab} or finite fields
\cite{Wie,LaM-Odl,PomS,Tei}.  The differences are because we need to
compute~$P$ and~$Q$, as explained in the introduction.  For matrices
over~$\Z$, one can find the SNF~$D$ matrix efficiently by reducing
modulo a number of primes \cite{DSV} \cite{DEvGU07}, or by other
techniques \cite{HHR93} \cite{Magma}.  Yet it is not clear how to
find~$P$ and~$Q$ by reducing by multiple primes.  The need for the
change of basis matrices is why \shh\ aims to work globally.

{\em Fill-in\/} is a concern in sparse linear algebra over any
ring~$R$.  Imagine two vectors that both have a nonzero at~$i_0$.  Add
a scalar multiple of the first to the second in order to zero out the
value at~$i_0$.  In general, the result will have nonzeros in the
union of the positions where the originals had nonzeros, apart
from~$i_0$.  (We follow the standard convention in the sparse matrix
literature and abbreviate ``nonzero entries'' to ``nonzeros.'') For
very sparse vectors, the number of nonzeros almost doubles.  Fill-in
is this growth in the number of nonzeros.

A separate issue when~$R$ is infinite is {\em integer explosion}.
Over~$\Z$, the length (number of digits) of the product of two numbers
is roughly the sum of the lengths of the factors.  A vector operation
that zeroes out one position will tend to increase the length of the
numbers in the rest of the vector.  \shh's purpose is to avoid
fill-in and integer explosion as much as possible.  With $R = \F_p$,
integer explosion is not an issue, and the focus is on avoiding
fill-in.

We find the SNF by performing a sequence of {\em elementary
operations\/}: permuting two columns, adding a scalar multiple of one
column to another, and multiplying a column by a unit of~$R$, plus the
same for rows.  The algorithm is described in~(\ref{SNF-pseudocode}).

Algorithms that operate on only one side of the matrix are more
familiar.  These include the Hermite normal form (HNF) $A = HQ$
\cite[(2.4.2)]{Cohen}.  Over a field, HNF is the same as Gaussian
elimination on columns, with $H$ in column-echelon form and
$Q\in\GL_n(\Z)$.

In principle, we prefer SNF to HNF because we are working with cochain
complexes.  To evaluate Hecke operators on~$H^i$, we need to compute
with the map $\mmker d^i \to (\mmker d^i)/(\im d^{i-1})$ that reduces
cocycles modulo coboundaries.  This requires the~$P$ matrix of
$d^{i-1}$ and the~$Q$ matrix of $d^{i}$.  When computing all of
$H^*$, it is natural to compute~$P$
and~$Q$ for~$d^*$ at the same time.

When the $d^i$ are very large, however, we must compromise
by omitting computation of the change-of-basis matrices that
are not needed.  Since this paper is about $H^5$, we compute for $d^5$
only $D$ and $Q$, and for $d^4$ only $P$ and $D$.  The biggest savings
arise because the largest matrices, $d^4$, are significantly wider
than their height, as Figure~\ref{fig-d4-d5} shows.  
The~$Q$ matrices for $d^{4}$, those on the longer side, are
fortunately the ones we can forget.

HNF does have one advantange over SNF when one is forgetting the $Q$
matrix: it can be computed by the following {\em disk HNF\/}
algorithm.  Write the whole matrix~$A$ to disk, then read it back in
one column at a time.  As one reads each column, one puts the matrix
accumulated so far into HNF.  Over $R = \F_p$, this means using
standard Gaussian elimination on columns, with no special pivoting
algorithm.  Again, the savings arise because $d^4$ has a short side
and a long side.  The echelon form never exceeds $n_5 \times n_5$, the
square on the short side.

With the wrong matrices, though, disk HNF is a recipe for disaster.
It can succeed only if the matrix has low co-rank.  The {\em
co-rank\/} of an $m\times n$ matrix is $\markmin(m,n)$ minus the
rank~$\rho$ of the matrix.  Assume $m \le n$ from now on (this holds
for our~$d^5$ and~$d^4$), so the co-rank is $m-\rho$.  Imagine that,
after some work, one must put such a matrix into column-echelon form
using Gaussian elimination.  We claim that the echelon form will have
a great deal of fill-in, no matter how cleverly the pivots are chosen.
The echelon form will have $\rho$ pivot rows with only a single
nonzero.  The remaining $m-\rho$ rows will in general be
dense---no pivot has had a chance to clear out their entries, and by
the law of averages they will mostly be nonzero.  Hence there are
about $(m-\rho) \cdot \rho = {}$ (co-rank)${}\cdot{}$(rank) nonzeros
in the echelon matrix.  We cannot stop~$\rho$ from being large.  But
when $m-\rho$ is also large, the product $(m-\rho) \cdot \rho$ is
prohibitively large.  These observations are for the final result of a
computation; usually fill-in is even worse in the middle of the
computation, before all the pivot rows have been established.

The main technical observation of this paper is to use the change of
basis matrices in a simple way to transform~$d^4$ into an equivalent
matrix~$\eta$ of low co-rank.  We start with SNF on~$\eta$, switch to
disk HNF when the fill-in of SNF forces us to do so, and rely on the
low co-rank to make disk HNF succeed.  The matrix~$\eta$ is defined in
Equation~\ref{def-eta} below.

In~(\ref{subsec-SNF}) and~(\ref{subsec-coh}) we present these ideas
more precisely.

\subsection{Computing the SNF}
\label{subsec-SNF}
Let $A$ by an $m\times n$ matrix with entries in a field~$R$ where
exact computation is possible.  We define {\sl elementary matrices}
$P_l \in \GL_m(R)$ as usual~\cite[p.~182]{Jac}.  These are
permutation, translation, and dilation matrices, corresponding to the
elementary operations listed above.  Replacing $A$ with $P_l A$ or
$P_l^{-1} A$ performs an elementary row operation on~$A$.  Multiplying
on the right by an elementary matrix $Q_l \in \GL_n(R)$ performs an
elementary column operation.

\subsubsection{The Markowitz algorithm}
\label{markowitz}
Markowitz \cite[(7.2)]{DER} \cite{HHR93} is a well-established approach to reducing
fill-in.  It is a greedy algorithm, reducing fill-in as much as
possible at each step.  Let $a_{ij}$ be a nonzero.  Let
$r_i$ be the number of nonzeros in $a_{ij}$'s row, and $c_j$
the number in its column.  If one adds a multiple of row~$i$ to row~$k$
in order to zero out the $a_{kj}$ entry, one creates up to $r_i - 1$
new nonzeros in row~$k$.  Using row~$i$ to clear out the entire
column~$j$ produces up to $(r_i - 1)(c_j - 1)$ new entries.  The
Markowitz algorithm, in its simplest form, chooses at each step the
pivot $a_{ij}$ that minimizes the {\em Markowitz count\/} $(r_i -
1)(c_j - 1)$.  (If~$R$ were~$\Z$, we would also need to address
integer explosion, by avoiding pivots with large absolute value even
if they have low Markowitz count.)

It can be slow to compute the Markowitz count for all $a_{ij}$.  One
approach is to look at only a few rows with small $r_i$---say the
smallest ten rows---and minimize the Markowitz count only for those
rows.  Early versions of the \shh\ code used this approach.
Currently, we prefer to avoid fill-in at whatever cost in speed, and
we always search over all entries.  To speed up the search, we store
the~$r_i$ and~$c_j$ in arrays and update the arrays with every
elementary operation.

\subsubsection{Statement of the algorithm}
\label{SNF-pseudocode}
We now describe \shh~2.2's SNF algorithm.  Implementation details are
deferred to~(\ref{implementation}).

The main strength of the algorithm is the interplay between the
Markowitz count and disk HNF when the co-rank is low.  Outside these
aspects, it is like many SNF algorithms \cite[(3.7)]{Jac}
\cite[(2.4.4)]{Cohen}.

The algorithm uses an index~$c$, the {\em pivot index\/} or {\em
  corner\/}, which just says the current pivot is at $a_{cc}$.  The
{\em active region\/} is where $c \le i < m$ and $c \le j < n$.
Outside the active region, the matrix has nonzeros on the diagonal
and nowhere else.

The parameter~$\tau$ controls when we switch from SNF to disk HNF.  It is chosen by the user based on heuristics and experience; see~(\ref{speed-space}) for details.
The part of the algorithm with~$\tau$ is stated when we are forgetting
the~$Q$'s, as we always do for~$d^4$; it would be easy to extend this
part for the~$P$'s also.

If one does not need to remember~$P$ and~$Q$, one simply
omits the part of the algorithm that writes them out.  Our
implementation overwrites $A$ with~$D$.

{\em Input:\/} an $m\times n$ sparse matrix~$A = (a_{ij})$ with
entries in a field~$R$.  If we are not finding the~$Q$ change of basis
matrices, we are given in addition a parameter~$\tau$ (defaulting to
$\tau = \infty$).

{\em Output:\/} A finite list $(P_0, P_1, P_2, \dots)$ of $m\times m$
elementary matrices, a finite list $(\dots, Q_2, Q_1, Q_0)$ of
$n\times n$ elementary matrices, and an $m\times n$ diagonal
matrix~$D$, such that
\begin{displaymath}
A = (P_0 P_1 P_2 \cdots ) \cdot D \cdot (\cdots Q_2 Q_1 Q_0)   \\
  = P D Q
\end{displaymath}
as in Equation~\ref{def-SNF}.

{\em Step 0.} Set $c = 0$.

{\em Step 1.}  If the number of nonzeros in the active
region is $\ge\tau$, and if disk HNF has not run yet, run disk HNF on
the active region.  That is, write the active region to disk, find
its column-echelon form while reading it back one column at a time,
and replace the old active region in memory with the echelon form, keeping track of the~$Q$'s.

{\em Step 2.}  If the active region is entirely zero, including the
case $c = \markmin(m,n)$, then return the lists of $P$'s and $Q$'s,
return~$A$ overwritten with~$D$, and terminate the algorithm.

{\em Step 3.} Choose a nonzero in the active region that
minimizes the Markowitz count.  This is the {\em pivot\/}.  Use a row
and column permutation if necessary to move the pivot to the $a_{cc}$
position (this is {\em complete pivoting\/}).  If the row permutation
is $A \to P_l^{-1}A$, then append $P_l$ to the right side of the list
of~$P$'s.  (Of course $P_l = P_l^{-1}$ for order-two permutations.)
Similarly, append the column permutation $Q_l$ to the left side of the
list of~$Q$'s.

{\em Step 4.}  For all $j$ with $c < j < n$ and $a_{cj} \ne 0$,
subtract a multiple of column~$c$ from column~$j$ to make $a_{cj} =
0$.  For each of these elementary operations $A \to A Q_l^{-1}$,
append $Q_l$ to the left side of the list of~$Q$'s.

{\em Step 5.}  For all~$i$ with $c < i < m$ and $a_{ic} \ne 0$,
subtract a multiple of row~$c$ from row~$i$ to make $a_{ic} = 0$.  For
each of these elementary operations $A \to P_l^{-1} A$, append $P_l$
to the right side of the list of $P$'s. (If~$R$ were not a field,
steps~4 and~5 would need to be extended when $a_{ij} / a_{cc}$ has a
nonzero remainder for some $a_{ij}$ in the active region.)

{\em Step 6.}  Increment~$c$ and go to step~1.

\subsubsection{Representing change-of-basis matrices}
\label{p-q-on-disk}
It is important that we return $P$ and $Q$ as lists of elementary
matrices.  The products $P_0 P_1 P_2 \cdots$ and $\cdots Q_2 Q_1 Q_0$
are likely to be dense; we could never hold them in RAM, much
less compute their inverses.  Fortunately, it is easy to work with
them as lists.  Given a matrix~$B$, compute $QB$ by multiplying~$B$ on
the left by~$Q_0$, $Q_1$, and so on.  To compute $Q^{-1} B$, run
through the~$Q_l$ in the opposite order, decreasing~$l$, and
multiply~$B$ on the left by $Q_l^{-1}$.  Similar comments apply to
the~$P_l$, and to transposes.

The lists are still too big to fit in RAM, so we store them on disk.
We read once through them every time we need to operate with~$P$
or~$Q$.  We use a text-based data format where each elementary matrix
takes up only about~20 characters.  Storing a translation matrix, for
example, only requires storing the two indices and the value of the
off-diagonal entry.  Reading the elementary matrices in left-to-right
order is straightforward.  To read in right-to-left order, we use a
pointer that steps backward through the file.

\subsubsection{Sizes of the matrices}
\label{n-i-estimates}
We will need more precise estimates for the $n_i = \dim \cochains^i$.
We refer the reader to~\cite{AGM1} for the definitions.  The~$n_i$ are
approximated by sums of terms of the form $|\P^3(\Z/N\Z)| /
|\Gamma_\sigma|$.  Here $\P^3 = \P^3(\Z/N\Z)$ is projective
three-space over the ring $\Z/N\Z$, and $\Gamma_\sigma$ are the
stabilizer groups of certain basis cells in the well-rounded
retract~$W$.

\begin{eqnarray*}
n_6 & \approx & \frac{1}{96} |\P^3| \\
n_5 & \approx & \left(\frac{1}{24} + \frac{1}{24} + \frac{1}{60}\right) |\P^3| = \frac{1}{10} |\P^3| \\
n_4 & \approx & \left(\frac{1}{8} + \frac{1}{6} + \frac{1}{24} + \frac{1}{72}\right) |\P^3| = \frac{25}{72} |\P^3|
\end{eqnarray*}

If~$N$ has the prime factorization $\prod p_i^{e_i}$, then
\begin{displaymath}
|\P^3| = N^3 \prod \frac{1 + p_i + p_i^2 + p_i^3}{p_i^3}.
\end{displaymath}
When~$N$ is prime, this reduces to the familiar formula $|\P^3| = 1 +
N + N^2 + N^3$.  In general, if we consider the set $\{p_i\}$ of
primes dividing~$N$ to be fixed, then $|\P^3|$ is a constant
times~$N^3$, where the constant depends on the set.  In the range
of~$N$ we consider, up to the low 200's, the largest constant we
encounter is 4.04, for $N = 2 \cdot 3 \cdot 5 \cdot 7 = 210$.  This
factor of~4 is why we said that our estimates for the~$n_i$ are a few
times larger for composite~$N$ than for prime~$N$.

The~$\approx$ symbols arise because the orbits of $\Gamma_\sigma$
in~$W$ are generically of size $|\Gamma_\sigma|$ but are occasionally
smaller.  What is more, we only count the orientable orbits.
Experiments with~$N$ in the range~40 to~100 suggest that the error in
the~$\approx$ is at worst 0.3\% for prime~$N$ and 0.6\% for
composite~$N$.

%%%   N    largest error in n_i for i=4,5,6
%%%   42   0.6%
%%%   48   0.6%
%%%   50   0.2%
%%%   55   0.3%
%%%   89   0.3%
%%%   97   0.3%

\subsection{Computing cohomology}
\label{subsec-coh}
Let $d^5$ and $d^4$ be as in Formula~\ref{main-complex}, with
$n_i = \dim \cochains^i$.  We now describe how we use \shh\ to compute
$H^5$ of the complex.

First, compute the SNF $d^5 = (?) D_5 Q_5$ with $P_5$ omitted.  Let
$\rho_5 = \rank D_5 = \rank d^5$.

Second, define
\begin{equation}
\eta = Q_5 d^4.  \label{def-eta}
\end{equation}
Since $d^*$ is a cochain complex, the topmost $\rho_5$ rows of~$\eta$
are zero.  Delete these rows, still calling the matrix~$\eta$.
Thus~$\eta$ has dimension $(n_5 - \rho_5) \times n_4$.

Third, compute the SNF $\eta = P_\eta D_\eta (?)$, with~$Q_\eta$
omitted.  Let $\rho_\eta = \rank D_\eta = \rank\eta$.  Note that
$\rank d^4 = \rho_\eta$, since~$d^4$ and~$\eta$ differ only by a
change of basis and deleting some zero rows.

We can now report the Betti number $h_5 = \dim H^5$:
\begin{displaymath}
h_5 = n_5 - \rho_5 - \rho_\eta.
\end{displaymath}
We need not only the Betti number, though, but an explicit list $z_1,
\dots, z_{h_5}$ of cocycles in $\mmker d^5$ that descend modulo~$\im
d^4$ to a basis of the cohomology.  Let $B$ be the $(n_5 - \rho_5)
\times h_5$ matrix with the identity in the bottom $h_5 \times h_5$
block and zeroes elsewhere.  Compute $\bar B = P_\eta B$.  Add to the
top of $\bar B$ the $\rho_5$ rows of zeroes that we deleted
from~$\eta$, still calling it~$\bar B$.  Then the columns of $Q_5^{-1}
\bar B$ are an appropriate list of cocycles~$z_j$.

Our Hecke operator algorithm takes a cocycle~$y$ and maps it to its
Hecke translate, a cocycle~$y'$.  For simplicity, assume $y = z_j$.
The Hecke translate~$z_j'$ is generally not in the span of the $z_j$.
Rather, $z_j' = s_{1,j} z_1 + \cdots + s_{h_5,j} z_{h_5} + b_j$, where
the $s_{i,j} \in R$ are scalars and~$b_j \in \im d^4$ is a coboundary.
Computing the~$s_{i,j}$ and~$b_j$ is straightforward, since $Q_5$,
$P_\eta$, and the $z_j$ are all stored on disk.  Ultimately, we
express each Hecke operator as the $h_5 \times h_5$ matrix $(s_{i,j})$
with respect to our cohomology basis.

\subsubsection{Co-rank of~$\eta$}
Using $\eta$ may seem inelegant because it violates the symmetry of a
cochain complex.  Since the complex is
\begin{equation*}
0 \longleftarrow \cochains^6
  \stackrel{P_5 D_5 Q_5}\longleftarrow \cochains^5
  \stackrel{P_4 D_4 Q_4}\longleftarrow \cochains^4
  \longleftarrow \cdots,
\end{equation*}
it is more elegant to compute $H^5$ using $Q_5$ and $P_4$, which are
both $n_5 \times n_5$ matrices.

However, $\eta$ has one great virtue: by removing the rows of zeroes
from its top, we have dropped its co-rank down to~$h_5$.  We observe
that $h_5$ is never more than~80 for the~$N$ we consider,
while~$\rho_5$ could reach into the millions.  This difference is what
allows disk HNF to succeed.

Let us give more precise estimates.  The Betti number $h_6 = \dim H^6$
equals $n_6 - \rho_5$ since our chain complex has only 0's after
degree~6.  We observe that~$h_6$ is never more than about~40 in our range
of~$N$.  Thus $\rho_5 \approx N^3/96 - 40$.  Estimating~$h_5$ as~80,
the rank $\rho_\eta = \rho_4 \approx N^3/10 - (N^3/96 - 40) - 80$.
Both~40 and~80 are negligible for large~$N$, so $\rho_\eta \approx
(1/10 - 1/96) N^3 = 43/480 N^3$.  For~$\eta$, the co-rank is again
about~80, meaning the number of entries in~$\eta$'s echelon form is
(co-rank)${}\cdot{}$(rank) ${} \approx 80 \cdot \rho_\eta = 80 \cdot
(43/480) N^3 \approx 7 N^3$.  But the number of entries in $d^4$'s
echelon form is $\approx (n_5 - \rho_4) \rho_4 \approx ((1/10)N^3 -
(43/480) N^3)\cdot((43/480)N^3) = (1/96)(43/480)N^6 \approx 0.0009
N^6$.  At $N \approx 200$, the latter is $\approx 7500 N^3$.  In other
words, the echelon form of~$d^4$ has about 1000 times more entries
than the echelon form of~$\eta$ when~$N$ is near~200.

%%% h_6 is 41 at level 48.  For prime N it is more like 10 or 20.  It seems
%%% to be biggest when N is highly divisible by 2.

\begin{figure}
\input{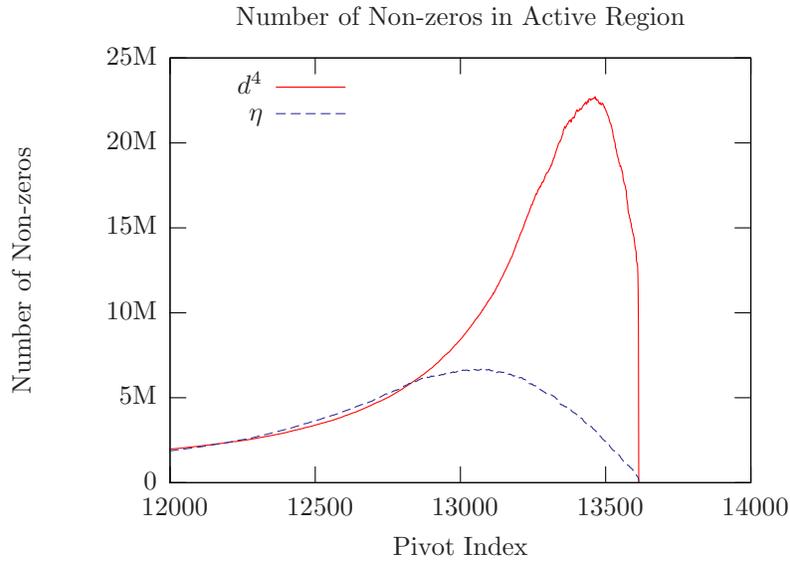}
\caption{Size of the active region during an SNF computation
for~$d^4$ and~$\eta$.  Here M denotes {million}.}
\label{fig-d4-n53}
\end{figure}

This analysis was for Gaussian elimination, not SNF.  When we compute
the SNF of matrices of large co-rank, we observe the same behavior
empirically.  Figure~\ref{fig-d4-n53} compares the fill-in for the SNF
computations of~$d^4$ and~$\eta$ at the same level $N = 53$.
Both matrices have 52766 columns and rank 13614.  The example uses
Markowitz only, not disk~HNF.  We show only the pivot indices $c \ge
12000$, since the graphs look the same to the left of that point.  The
fill-in for~$d^4$ is clearly the worse of the two, with a peak over
three times higher than for~$\eta$.  In general, the SNF algorithm
displays ``catastrophic cancellation'': the number of nonzeros in
the active region tends to grow rapidly
%\footnote{[As
%Figure~\ref{fig-d4-n53} shows, $d^4$ grows slowly, then turns a sharp
%corner.  The corner is too sharp for the overall curve to be either
%exponential or quadratic.  Maybe it's linear while the columns are so
%sparse that each pivot adds $(6-1)(26-1) = ({\rm constant})$ rows,
%linear like this for a long time, then segues into exponential.  See
%if \cite{PomS} has any ideas.]} 
until almost the end of the algorithm, when the number decreases
sharply to zero.  Catastrophic cancellation begins for~$d^4$ at
row~13464 and for~$\eta$ at about row~13084.

\begin{figure}
\input{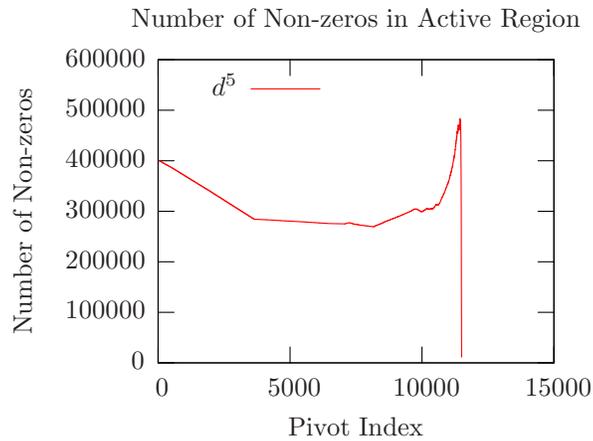}
\caption{Size of the active region during an SNF computation for~$d^5$.}
\label{fig-d5-n103}
\end{figure}

The fill-in for our smaller matrix~$d^5$ is harder to predict.  There
are many columns with only one or two entries.  These allow Markowitz
to reduce the matrix with no fill-in at all, at least in the initial
stages.  Later, the number of nonzeros grows rapidly as for~$d^4$,
with an even more precipitous cancellation at the end.
Figure~\ref{fig-d5-n103} gives the example for level $N = 103$.

\subsection{Implementation} \label{implementation} \shh~2.2 is written
in Common Lisp, which is the ideal language for it in many ways.
Common Lisp source code compiles to very efficient machine code;
carefully-written Lisp is as fast as C++~\cite{cll}.  Yet writing
high-level code in Lisp is easier than in C or C++.  Like its
inheritor Python, Lisp has an outer real-eval-print loop that promotes
high-level thinking and rapid development.  Type safety is built in.
Lisp's object-oriented package CLOS is the most flexible of any
language's.  Arbitrary-precision integers and rationals are part of
the language standard and can be very fast.

\shh~2.1 and~2.2 were developed with Allegro Common Lisp (ACL) from
Franz, Inc.  Under Linux we use the free implementation Steel Bank
Common Lisp (SBCL).  \shh~2.0 was developed in Java.  \shh~1.$x$ was
developed in Carnegie-Mellon University Common Lisp (CMUCL), the
parent of SBCL.  All the Lisps mentioned here produce high-quality
compiled code.

\shh~1.$x$ was restricted to $R = \Q$, but had general support for the
derived category of complexes of sheaves over cell complexes.  It was
also ported to Macaulay~2~\cite{M2}.

\subsubsection{Overview of data structures} \label{impl-overview} Any
sparse matrix package uses data structures that store only the
nonzeros.  \shh\ stores a sparse matrix as an array of sparse vectors
representing the columns.  A sparse vector is a linked list of
nonzero sparse elements.  A {\em sparse element\/} represents an
ordered pair $(i,v)$ with $i\in\N$, $v\in R$.  When this sparse
element occurs in column~$j$, it means the $i,j$ entry of the matrix
has nonzero value~$v$.  The rest of~(\ref{implementation}) gives
further details.

\subsubsection{Testing platforms}
To explain the test data we will present later, we give some details
about the two machines we use most.  Portland, our Windows machine, is
a laptop with a 1.30~GHz Celeron~M processor and 992~MB of RAM,
running ACL~8.1 under Windows XP 2002.  Athens, our Linux machine, is
a 64-bit laptop with a 2.20~GHz Intel Core2 Duo processor and 3.9~GB
of RAM, running SBCL~1.0.12 under Linux 2.6 (Ubuntu~8).

\subsubsection{Implementation of sparse elements}
\label{sparse-elt}
For general rings~$R$, \shh\ implements the indices $i\in\N$ and
values $v\in R$ as objects.  They are allocated on the heap, and
contain type information needed to use them in an object-oriented way.
For these rings, we store $(i,v)$ as a pair of pointers to~$i$ and
to~$v$.  The pair of pointers is a native data structure in Lisp, the
{\em cons\/}, the fundamental building block for linked lists and tree
structures.

The cons implementation is convenient, but, as in all languages, the
indirection carries a penalty in efficiency.  Lisp may store the
numbers on a different memory page from the conses themselves, forcing
the machine to flip back and forth between pages.

When $R=\F_p$, we implement sparse elements more efficiently.  We
typically use primes~$p$ a little less than~$2^{15}$, so that sums and
products in $[0,p)$ can be computed in a 32-bit integer before
reducing mod~$p$.  Let~$k$ be the smallest integer such that $p <
2^k$.  (For our favorite prime 12379, $k = 14$.)  Each $v\in\F_p$ fits
in~$k$ bits.  We would like to pack $(i,v)$ into a single 32-bit
integer, but we cannot.  With $k=14$, say, there would only be 18 bits
left for the index, and we need a larger range of indices than $2^{18}
\approx 10^6$.  Therefore, on a 32-bit machine, we pack $(i,v)$ into
an arbitrary-precision integer, but arrange things so the integer will
fit into 32 bits at the critical stages of the computation.  Let~$M$
be an upper bound on the number of rows in the matrix.  We store
$(i,v)$ as the integer
\begin{equation}
(M-i-1)\cdot 2^k + v.  \label{def-12379x}
\end{equation}
Near the end of a computation, when space is critical, $i$~will be
close to~$M$.  Hence $M-i-1$ will be small and
Formula~\ref{def-12379x} will fit in 32 bits.

On a 64-bit machine, we store $(i,v)$ even more simply as
\begin{equation}
i\cdot 2^k + v.  \label{def-12379-64}
\end{equation}
For us, this never exceeds a 64-bit integer.  When $k = 14$, for
instance, it would exceed it only if $i > 2^{64-14} = 2^{50} \approx
10^{15}$, whereas our largest~$m$ are around~$10^6$ or~$10^7$.  By
declaring the sparse element to be a 64-bit integer throughout the
program, we avoid the separate heap allocations for~$i$ and~$v$, the
pair of pointers holding~$i$ and~$v$, and the object-oriented
indirection.

In all of \shh's implementations of sparse elements, operations in the
ring~$R$ really operate on the sparse elements.  For instance, the sum
of $(i_1, v_1)$ and $(i_2, v_2)$ is a new sparse element with value
$v_1 + v_2 \in R$ and with some~$i$ determined by the context,
typically~$i = i_1$.

\subsubsection{Implementation of sparse vectors and matrices}
\shh's sparse matrix data structure is an array of sparse vectors
giving the columns, together with arrays of the row lengths~$r_i$ and
column lengths~$c_j$ used for Markowitz.  \shh\ implements a sparse
vector as a singly linked list of sparse elements, sorted by
increasing index~$i$.  The backbone of the linked list (as opposed to
the sparse elements it holds) is again built from conses.
Implementing vectors as a linked list is flexible, compared to
constantly allocating and reallocating arrays in~C.  Lists grow in
size as fill-in occurs, and shrink as entries cancel each other out.
As they shrink, the memory is reclaimed by the garbage collector.

\shh\ includes a macro \texttt{with-splicer} as a mini-language for
surgery on lists.  \texttt{with-splicer} iterates down a list and offers
read, insert, modify and delete operations for individual elements,
plus splicing commands to add and remove larger chunks.  We allocate
as few new backbone conses as possible.

Implementing sparse vectors as singly-linked lists is flexible, as we
have said, but it involves the usual risks because accessing a single
element might take time linear in the length of the list.  We can
avoid this trap with a little care.  To find the sum {\bf x} + {\bf y}
of sparse vectors {\bf x} and {\bf y}, for example, we run through the
vectors in parallel and accumulate the sum as we go.  Similar comments
apply to scalar multiplication, subtraction, and the dot product.

One place we might encounter quadratic behavior is step~5 of the SNF
algorithm~(\ref{SNF-pseudocode}).  If the pivot row and column have
$r_i$ and $c_j$ entries respectively, a direct approach would require
$(r_i - 1)(c_j - 1)$ operations, each linear in the length of the
column.  The trick here is simply to put step~4 before step~5.  Step~4
handles the column operations linearly as in the previous paragraph,
and step~5 then has $r_i = 1$.

Another likely place for slow behavior is multiplying a sparse matrix~$B$ on
the left by an elementary matrix.  This includes the important
computation of~$\eta$ (Equation~\ref{def-eta}).  An elementary
operation with rows~$i$ and~$j$ might involve a linear sweep down $r_i
+ r_j$ columns.  We handle this situation by taking the transpose
of~$B$, multiplying on the right by the transposes of the elementary
matrices, and taking a final transpose to put everything back at the
end.  The transpose of~$B$ takes time linear in the number of entries,
which for a sparse matrix is small by definition.

\subsubsection{Comparison of sparse element implementations}
\label{sparse-elt-compare} Table \ref{speedtable} shows a test of
speed and space requirements for the three implementations of sparse
elements over~$\F_{12379}$ on our Linux machine.  We timed the SNF
computation for the~$d^4$ matrix for level $N = 53$, the matrix of
Figure~\ref{fig-d4-n53}.  The matrix is $15218 \times 52766$.  We used
only the Markowitz portion of~(\ref{SNF-pseudocode}), no disk HNF.
Since we were on a 64-bit machine, each sparse element in
Formula~(\ref{def-12379-64}) takes 8 bytes.
Formula~(\ref{def-12379x}) has about the same space requirement,
especially towards the end of a computation when~$i$ is close to~$M$.
The $(i,v)$ representation requires 8 bytes for each of~$i$ and~$v$,
plus 16 bytes for the cons itself, for a total of 32 bytes.  To all
these figures we must add 16 bytes for the cons in the backbone of the
sparse vector that holds the sparse element.

\subsubsection{Overall speed and space requirements}
\label{speed-space}
To summarize, our implementations of sparse elements are optimized for
both space and speed, and our sparse vector and matrix algorithms
avoid speed traps leading to quadratic time unnecessarily.  On the
other hand, at the higher layers of the algorithm, we sacrifice speed,
minimizing fill-in at all costs.  For instance, we mentioned
in~(\ref{markowitz}) that we do a full Markowitz scan at each pivot
step.  This takes about one third of the total time for the algorithm
until we switch to disk~HNF.

The largest matrix we have handled so far has size $845,712 \times
3,277,686$.  This is the~$\eta$ matrix for level $N = 211$.  It has
close to 20~million nonzeros.  We carried out the computation on our
Linux machine, using the implementation of~$\F_{12379}$ in
Formula~\ref{def-12379-64}.  The sizes of $d^5$ and $d^4$ are $98,351
\times 944,046$ and $944,046 \times 3,277,686$, respectively.  We
reduced~$d^5$ using only the Markowitz portion
of~(\ref{SNF-pseudocode}), with no disk HNF.  We reduced~$\eta$ using
both Markowitz and disk HNF, switching when the active region had
$\tau$ equal to about 116~million nonzeros.  Converting~$d^4$
to~$\eta$ as in~(\ref{def-eta}) took significant time by itself, since
it called for over three million elementary column operations on the
transpose of~$d^4$.  How the running time broke down is shown in Table
\ref{runningtable}.

%%% Here is how I computed the timing for N = 211.  d^5 and d^4 were
%%% written out as delta?bak.txt at 6:21 on Sep 26.  The reduction of d^5,
%%% which was 98351 by 944046, finished at 12:57 on the same day.  It
%%% spent an amount of time TBD converting d^4 to eta.  Reduction of
%%% eta, with disk-hnf-threshold 0.0005, began at 2:24 on day X.  Day X+1
%%% at 2:24 was about corner 103600.  Day X+2 (all until further notice
%%% are at 2:24) was about corner 190500.  Day X+3 was 269300.  Day X+4
%%% was 341800.  Day X+5 was 406500.  Day X+6 was 465700.  Day X+7 was
%%% 530000.  Day X+8 was 587700.  Day X+9 was 641400.  Day X+10 was
%%% 687400.  Day X+11 was 722800.  Day X+12 was 749450.  At 6:32 on Day
%%% X+12, it wrote out the Markowitz part at corner 753854 and handed the
%%% work off to disk HNF.  It took a further day, till 6:43 on Day X+13,
%%% to reverse the file of P matrices.  This was on Oct 15, so Day X was
%%% Oct 2.  The time to convert d^4 to eta (look up how many elementary
%%% operations it was) was thus from 12:57 on Sep 26 to 2:24 on Oct 2,
%%% which is 5~1/2 days.  Disk HNF took about 13 more days.

It is interesting to compare our running times with those
in~\cite{DEvGU07}.  They compute for $\GL_7(\Z)$ at level $N = 1$,
while we do $\SL_4(\Z)$ at $N = 211$.  The number of nonzeros
is roughly comparable, 37~million versus 20~million.  Both
computations took about one month.  They computed mod~$p$ for
several~$p$, but used 50~processors; we did one~$p$ on one processor.
We find ourselves joking that
\[
\GL_7 = \GL_4 + 211^3.
\]
How the running time broke down is shown in Table \ref{runningtable}.
We do not distinguish between the wall-clock time and CPU time because
they are essentially the same.  We ran the computation on our Linux
machine on a single processor.  The machine was doing nothing
besides the computation.  Occasional tests with {\texttt{top}} showed the
CPU running consistently at 100\%.  We presume one of the two cores
ran the computation, and the other took care of background jobs like
displaying the desktop.

Recall that $\tau$ is the maximum number of nonzeros allowed in the
active region before switching from Markowitz to disk HNF.  Table
\ref{tautable} shows the largest~$\tau$ we have used successfully.
They depend on the chosen implementation of sparse elements, as well
as on the operating system and version of Lisp.  A~+ means we have
relatively little data for this combination of parameters, and~$\tau$
could likely go higher than shown.  Values without a~+ represent a
reasonable maximum, determined by long experience and many
out-of-memory crashes.  The number of bytes is computed as
in~(\ref{sparse-elt-compare}) for our Linux machine.  On our Windows
machine, a 32-bit integer counts 4 bytes, a cons 8 bytes.

\begin{table}[ht]
\begin{center}
\begin{tabular}{|l||l|l|}
\hline
{\em implementation\/} & {\em total time\/} & {\em space per entry\/}
\\
\hline
$(i,v)$ as a cons & 2768 sec & 48 bytes \\
Formula~\ref{def-12379x} & 1385 sec & 24 bytes \\
Formula~\ref{def-12379-64} & 784 sec & 24 bytes \\
\hline
\end{tabular}
\end{center}
\medskip
\caption{\label{speedtable}Comparison of sparse element implementations.}
\end{table}

\begin{table}[ht]
\begin{center}
\begin{tabular}{|l|l|}
\hline
SNF of $d^5$ & $\frac{1}{4}$ day \\
converting $d^4$ to $\eta$ & $5\frac{1}{2}$ days \\
SNF of $\eta$, Markowitz portion & 12 days \\
SNF of $\eta$, disk HNF portion & 13 days\\
\hline
\end{tabular}
\end{center}
\medskip
\caption{\label{runningtable}Overall speed for level $N = 211$.}
\end{table}

\begin{table}[ht]
\begin{center}
\begin{tabular}{|l||l|l|l|l|}
\hline
{\em machine and RAM\/} & {\em implementation\/} & {\em
largest\/}~$\tau$ & {\em space per entry\/} & {\em total space\/} \\
\hline
Windows 1GB & $(i,v)$ as a cons & 30M & 24 & 720MB \\
Windows 1GB & (\ref{def-12379x}) & 22M+ & 12 & 264MB+ \\
Linux 4GB & (\ref{def-12379x}) & 42M+ & 24 & 1.0GB+ \\
Linux 4GB & (\ref{def-12379-64}) & 148M & 24 & 3.55GB\\
\hline
\end{tabular}
\end{center}
\medskip
\caption{\label{tautable}Maximum number of nonzeros allowed in the active region.}
\end{table}

\subsubsection{Other approaches}
\label{other-approaches}
We mention a few more sparse matrix techniques that appear in the
literature.

Many scientific applications involve sparse matrices with a pattern,
such as tridiagonal or banded diagonal.  The matrices in this paper
have no recognizable sparsity pattern.  A matrix coming from a
$d$-dimensional topological space would have a pattern in~$d$
dimensions, but not when flattened into a two-dimensional table of
rows and columns.

The RSA challenge matrices described in~\cite{PomS} had some special
properties.  The columns on the left were very sparse, and could be
used to clean out the somewhat denser columns on the right.  Rows and
columns with only two entries gave an even quicker reduction
step~\cite{LaM-Odl}.  The \cite{DEvGU07} matrices had many rows with
only one entry \cite[2.6.4]{DEvGU07}, a result of cancellation at
the low level $N=1$.  By and large, our matrices do not have these
properties.  The sparsity is almost entirely uniform.  The~$d^5$ have
a substantial fraction of columns with one or two entries, but
not~$d^4$.

Block-triangularization is another well-established technique for
sparse matrices.  Given an $m\times n$ matrix~$A$, we look for
permutation matrices $P_b$ and $Q_b$ so that $B = P_b A Q_b$ is {\em
block upper-triangular\/}: it has square blocks down the diagonal and
nothing but zeroes below the blocks.  The matrix can then be reduced
one block at a time, either to HNF or SNF.  The upper-triangular part
can be handled directly by back-solving.  Since we only permute the
matrix entries, there is no fill-in and no integer explosion.  Assume
for the moment that~$A$ is square and of full rank, and that after a
row permutation the diagonal~$a_{ii}$ is all nonzero.  For such~$A$,
the block-triangular form is unique, and the way to find it is well
known~\cite[Ch.~6]{DER}.  When~$A$ is not of full rank, the block
decomposition is generalized to the Dulmage-Mendelsohn decomposition,
which is roughly upper-triangular~\cite[(8.4)]{Davis}.  In our case,
$A$ is a fraction of a percent away from full rank and from having
nonzeros down the diagonal; for square~$A$ of this type, finding the
Dulmage-Mendelsohn decomposition takes about the same time and space
resources as block decomposition.  So far these algorithms are
polynomial time.  A new idea is needed, however, when~$A$ is not
square but rectangular, as it is for us.  One can find the
Dulmage-Mendelsohn decomposition of the left-hand $m\times m$ square,
then swap in columns from the wide section on the right in a way that
shrinks the left-hand diagonal blocks.  Deciding which columns to swap
is in general an NP-hard problem.  The third author has found good
decompositions of some rectangular~$A$ using a heuristic for which
columns to swap in.  One iterates the procedure ``do
Dulmage-Mendelsohn, then swap in columns'' many times.  We defer the
details to a future paper.

\section{Eisenstein cohomology and Paramodular forms}\label{s:eisenstein.paramodular}

In this section we provide the necessary background to state
Conjectures \ref{conj:eisenstein} and \ref{conj:siegel} and explain
our computational results in Section
\ref{s:conjectures.computational}.

\subsection{Hecke eigenclasses and Galois representations} 
We will describe some classes appearing in $H^{5} (\Gamma_{0} (N);
\C)$ in terms of the Galois representations conjecturally attached to
them.  Thus we begin by recalling what this means \cite{AGM2}.

Let $\xi \in H^{5} (\Gamma_{0} (N); \C)$ be a Hecke eigenclass.  In
other words, $\xi$ is an eigenvector for certain operators
\[
T (l,k) \colon H^{5} (\Gamma_{0}
(N); \C) \rightarrow H^{5} (\Gamma_{0} (N); \C),
\]
where $k=1,2,3$ and $l$ is a prime not dividing~$N$.  These operators
correspond to the double cosets $\Gamma_{0} (N) D (l,k) \Gamma_{0}
(N)$, where $D (l,k)$ is the diagonal matrix with $4-k$ ones followed
by $k$ $l$'s.  (One can also define analogues of the operators $U_{l}$
for $l\mid N$, although we do not consider them in this paper.)
Suppose the eigenvalue of $T (l,k)$ on $\xi$ is $a(l,k)$.  We define
the Hecke polynomial $H (\xi)$ of $\xi$ by
\begin{equation}\label{eqn:glheckpol}
H (\xi) = \sum_{k} (-1)^{k}l^{k (k-1)/2} a (l,k) T^{k} \in \C [T].
\end{equation}

Now we consider the Galois side.  Let $\Gal (\bar \Q / \Q)$ be the
absolute Galois group of $\Q$.  Let $\rho \colon \Gal(\bar \Q /
\Q)\rightarrow \GL_{n} (\bar\Q_{p})$ be a continuous semisimple Galois
representation unramified outside $pN$.  Fix an isomorphism $\varphi
\colon \C \rightarrow \bar\Q_{p}$. Then we say the eigenclass
$\xi$ is \emph{attached} to $\rho$ if for all $l$ not dividing $pN$ we
have
\[
\varphi (H (\xi))= \det (1-\rho (\Frob_{l})T),
\]
where $\Frob_{l}\subset \Gal (\bar \Q /\Q )$ is the Frobenius
conjugacy class over $l$.  Let $\varepsilon$ be the $p$-adic
cyclotomic character, so that $\varepsilon (\Frob _l) = l$ for any
prime $l$ coprime to $p$.  We denote the trivial representation by
$i$.

\subsection{Eisenstein cohomology}
Let $X = \SL_4 (\R)/ \SO (4)$ be the global symmetric space, and let
$X^{\BorelSerre}$ be the partial compactification constructed by Borel
and Serre \cite{bs}.  The quotient $Y := \Gamma _{0} (N)\backslash X$
is an orbifold, and the quotient $Y^{\BorelSerre } := \Gamma_{0}
(N)\backslash X^{\BorelSerre}$ is a compact orbifold with corners.  We
have 
\[
H^{*} (\Gamma_{0} (N); \C) \simeq H^{*} (Y; \C) \simeq H^{*}
(Y^{\BorelSerre}; \C).
\]

Let $\partial Y^{\BorelSerre} = Y^{\BorelSerre}\smallsetminus Y$.  The
Hecke operators act on the cohomology of the boundary $H^{*} (\partial
Y^{\BorelSerre}; \C)$, and the inclusion of the boundary $\iota \colon
\partial Y^{\BorelSerre} \rightarrow Y^{\BorelSerre}$ induces a map on
cohomology $\iota^{*}\colon H^{*} (Y^{\BorelSerre}; \C) \rightarrow
H^{*} (\partial Y^{\BorelSerre}; \C)$ compatible with the Hecke
action.  The kernel $H^{*}_{!}  (Y^{\BorelSerre}; \C)$ of $\iota^{*}$
is called the \emph{interior cohomology}; it equals the image of the
cohomology with compact supports.  The goal of Eisenstein cohomology
is to use Eisenstein series and cohomology classes on the boundary to
construct a Hecke-equivariant section $s\colon H^{*} (\partial
Y^{\BorelSerre}; \C) \rightarrow H^{*} (Y^{\BorelSerre}; \C)$ mapping
onto a complement $H^{*}_{\Eis} (Y^{\BorelSerre}; \C )$ of the
interior cohomology in the full cohomology.  We call classes in the
image of $s$ \emph{Eisenstein classes}.  (In general, residues of
Eisenstein series can give interior, noncuspidal cohomology classes,
with infinity type a Speh representation, but as noted in \cite{AGM1},
these do not contribute to degree~5.)

\subsection{Paramodular forms}\label{ss:paramodular} We now give some
background on Siegel modular forms.  We will skip the basic
definitions, which can be found in many places (e.g.~\cite{vdg}), and
instead emphasize paramodular forms.  Since our main application is
cohomology of subgroups of $\SL_{4} (\Z)$, we will focus on $\Sp_{4}$.

Let $K (N)$ be the subgroup of $\Sp_{4} (\Q)$ consisting of all
matrices of the form 
\[
\left(\begin{array}{cccc}
\Z & N\Z & \Z & \Z \\
\Z & \Z &\Z &N^{-1}\Z \\
\Z &N\Z &\Z&\Z \\
N\Z &N\Z&N\Z &\Z   
\end{array} \right).
\]
The group $K (N)$ is called the \emph{paramodular group}.  It contains
as a subgroup the standard congruence subgroup $\Gamma '_{0} (N)$
modeled on the Klingen parabolic; that is, $\Gamma '_{0} (N) \subset
\Sp_{4} (\Z)$ is the intersection $K (N)\cap \Sp_{4} (\Z)$.
\emph{Paramodular forms} are Siegel modular forms that are modular with
respect to $K (N)$.  Clearly such forms are also modular with respect to 
$\Gamma '_{0} (N)$, although modular forms on the latter
are not necessarily paramodular.  Note also that in the embedding
$i\colon \Sp_{4}(\Z) \rightarrow \SL _{4} (\Z)$, we have $i (\Gamma
'_{0} (N)) = i (\Sp_{4} (\Z)) \cap \SL_{4} (\Z)$.

The paramodular forms of interest to us are those of prime level~$N$
and weight $3$.  We denote the complex vector space of such forms by
$P_{3}(N)$.  One can show that $P_{3} (N)$ consists entirely of
cuspforms, i.e.~there are no weight $3$ paramodular Eisenstein
series.  Recently T.~Ibukiayama \cite{ibuk2} proved a general
dimension formula for $P_{3} (N)$:

\begin{theorem}\label{thm:ibuk}\emph{\cite[Theorem 2.1]{ibuk2}}
Let~$N$ be prime and let $\kappa(a)$ be the Kronecker symbol
$(\frac{a}{N})$.  Define functions $f, g\colon \Z \rightarrow \Q$ by 
\[
f (N) = \begin{cases}
2/5&\text{if $N\equiv 2, 3 \bmod 5$,}\\
1/5&\text{if $N=5$,}\\
0&\text{otherwise},
\end{cases}
\]
and 
\[
g (N) = \begin{cases}
1/6&\text{if $N\equiv 5\bmod 12$,}\\
0&\text{otherwise}.
\end{cases}
\]
We have $\dim P_{3} (2) = \dim P_{3} (3) = 0$.  For
$N\geq 5$, we have 
\begin{align*}
\dim P_{3} (N) &= (N^{2}-1)/2880\\
&+ (N+1) (1-\kappa (-1))/64 
+ 5 (N-1) (1+\kappa (-1))/192 \\
&+ (N+1) (1-\kappa (-3))/72 + (N-1)
(1+\kappa (-3))/36 \\
&+ (1-\kappa (2))/8+f (N)+g (N) -1.
\end{align*}
\end{theorem}

For any~$N$, the space of weight $k$ paramodular forms contains a
distinguished subspace $P_{3}^{\g} (N)$ originally constructed by
Gritsenko \cite{grit}.  This space is defined by a lift from the
space $J^{\cusp}_{k,N}$ of cuspidal \emph{Jacobi forms} of
weight $k$ and index~$N$ to $P_{3} (N)$.  We will not define Jacobi
forms here, and instead refer the reader to \cite{ez} for background.
For our purposes, we will only need to know the dimension $\dim
P_{3}^{\g} (N) = \dim J^{\cusp}_{3,N}$.  Formulas for the dimensions
of spaces of Jacobi forms can be found in \cite[pp.~121, 131-132]{ez};
the following reformulation is due to N.~Skoruppa:

\begin{theorem}\label{thm:dimgrit}
We have 
\begin{equation}\label{eq:dimgrit}
\dim J^{\cusp}_{3,N} = \sum_{j=1}^{m-1} s(k+2j-1)-\Bigl\lfloor
\frac{j^{2}}{4m} \Bigr\rfloor,
\end{equation}
where $s (k)$ is the dimension of the space of cuspidal elliptic
modular forms of full level and weight $k$.
\end{theorem}

Let $P^{\ng}_{3} (N)$ be the Hecke complement in $P_{3} (N)$ of the
subspace $P^{\g}_{3} (N)$ of Gritsenko lifts.  The dimension of this
space is easily determined by Theorems~\ref{thm:ibuk} and
\ref{thm:dimgrit}.

We conclude our discussion of paramodular forms by defining the Hecke
polynomial attached to an eigenform.  More details can be found in
\cite{py2}.  Let $l$ be a prime not dividing~$N$.  Then associated
to $l$ there are two Hecke operators $T_{l}$ and $T_{l^{2}}$.  For an
eigenform $h\in P_{3} (N)$ we denote the corresponding eigenvalues by
$\delta_{l}$, $\delta_{l^{2}}$:
\[
T_{l}h = \delta_{l}h, \quad T_{l^{2}}h = \delta_{l^{2}}h.
\]
We define the Hecke polynomial attached to $h$ by 
\begin{equation}\label{eqn:sympheckpol}
H_{\Sp}  (h) = 1-\delta_{l}T+ (\delta_{l}^{2}-\delta_{l^{2}}-l^{2})T^{2} -
\delta_{l}l^{3}T^{3}+l^{6}T^{4}.
\end{equation}
This polynomial is essentially the local factor at $l$ attached to the
spinor $L$-function for $h$.

\section{Conjectures and computational
results}\label{s:conjectures.computational}

In this section we present two conjectures on the structure of $H^{5}
(\Gamma_{0} (N); \C)$ for~$N$ prime, and conclude by describing our
computational evidence for them. 

\subsection{Notation} We begin by fixing notation for the different
constituents of the cohomology.

\begin{itemize}
\item \emph{Weight two holomorphic modular forms:} Let $\sigma_{2}$ be
the Galois representation attached to a holomorphic weight 2 newform
$f$ of level~$N$ with trivial Nebentypus.  Let $\alpha$ be the
eigenvalue of the classical Hecke operator $T_{l}$ on $f$.  Let $\IIa
(\sigma_{2})$ and $\IIb (\sigma_{2})$ be the Galois representations in
the first two rows of Table \ref{grtab} (see p.~\pageref{grtab}).

\item \emph{Weight four holomorphic modular forms:} Let $\sigma_{4}$
be the Galois representation attached to a holomorphic weight 4
newform $f$ of level~$N$ with trivial Nebentypus.  Let $\beta$ be the
eigenvalue of the classical Hecke operator $T_{l}$ on $f$.  Let
$\IV(\sigma_{4})$ be the Galois representation in the third row of
Table \ref{grtab}.

\item \emph{Cuspidal cohomology classes from subgroups of $\SL_{3}
(\Z)$:} Let $\tau$ be the Galois representation conjecturally
attached to a pair of nonselfdual cuspidal cohomology classes
$\eta,\eta ' \in H^{3} (\Gamma^{*}_{0} (N); \C)$, where $\Gamma^{*}_0
(N)\subset \SL_{3} (\Z)$ is the congruence subgroup with bottom row
congruent to $(0,0,*)$ modulo~$N$.  Let $\gamma$ be the eigenvalue of
the Hecke operator $T_{l,1}$ on $\eta $, and let $\gamma '$ be its
complex conjugate.  Let $\IIIa(\tau)$ and $\IIIb (\tau)$ be the Galois
representations in the last two rows of Table \ref{grtab}.

%\item \emph{Weight three paramodular forms:} Let $\rho_{3}$ be the
%Galois representation attached to a weight 3 paramodular eigenform $h$ of level
%$N$.  Let $\delta_{l}$ and $\delta_{l^2}$ be the eigenvalues of the Hecke
%operators $T_{l}$ and $T_{l^{2}}$ as described in Section
%\ref{ss:paramodular}. 
\end{itemize}

If $f$ is a weight $2$ or weight $4$ eigenform as above, or a weight
$3$ paramodular eigenform, we denote by $d_{f}$ the degree of the
extension of $\Q$ generated by the Hecke eigenvalues of $f$.  We say
that two eigenforms $f, g$ are \emph{Galois conjugate} if there is an
automorphism $\sigma \in \Gal (\bar \Q /\Q)$ such that the Hecke
eigenvalues of $f$ are taken into those of $g$ by $\sigma$.  We say
$f, g$ are \emph{equivalent} if $g$ is a $\Q$-linear combination of
$f$ and its Galois conjugates.  We extend these notions in the obvious
way to eigenclasses $\eta \in H^{3} (\Gamma^{*}_{0} (N); \C)$.

For any modular form of weight $k$ with Fourier expansion $f (z) =
\sum_{n} a_{n}e^{2\pi i nz}$, let $L(s,f)$ be the Dirichlet series
$\sum_{n} a_{n}/n^{s}$.  The series $L (s,f)$ can be completed to a
function $\Lambda (s,f)$ satisfying a functional equation of the shape
$s \rightarrow k-s$.

\subsection{Eisenstein cohomology}

\begin{conjecture}\label{conj:eisenstein}
Let~$N$ be prime.  Then the cohomology group $H^{5} (\Gamma_{0} (N);
\C)$ contains the following Eisenstein subspaces:
\begin{enumerate}
\item For each equivalence class of weight two holomorphic newforms of
level~$N$, choose a representative $f$ with associated Galois
representation $\sigma_{2}$.  Then there are two $d_{f}$-dimensional
subspaces in the cohomology, one attached to the Galois representation
$\IIa (\sigma_{2})$, and the other to the Galois representation $\IIb
(\sigma_{2})$.
\item For each equivalence class of weight four holomorphic newforms of
level~$N$, choose a representative $f$ with associated Galois
representation $\sigma_{4}$.  Then if the central special
value $\Lambda  (2,f)$ vanishes, there is a $d_{f}$-dimensional
subspace in the cohomology attached to
the Galois representation $\IV (\sigma_{4})$.
\item For each equivalence class of nonselfdual cuspidal cohomology
classes in $H^{3} (\Gamma^{*}_{0} (p); \C)$, $\Gamma^{*}_0 (p)\subset
\SL_{3} (\Z)$, choose a representative $\eta$ and let $\tau $ be the
conjecturally associated Galois representation.  Then there are two
$d_{\eta }$-dimensional subspaces, one attached to the Galois
representation $\IIIa (\tau)$, and the other to the Galois
representation $\IIIb (\tau)$.
\end{enumerate}
Furthermore, for~$N$ prime this is a complete description of the
Eisenstein subspace of $H^{5} (\Gamma_{0} (N); \C)$.
\end{conjecture}

In our earlier paper \cite{AGM2}, we also gave a conjecture about some
Eisenstein subspaces of $H^{5}$.  In fact, for weight $2$ modular
forms and for $\SL_{3}$-cuspidal cohomology, there is no difference
between \cite[Conjecture 1]{AGM2} and the conjecture here.  The new
part is in the contribution of the weight $4$ modular forms.  In
\cite{AGM2}, our data was only sufficient to suggest that the weight
$4$ forms $f$ appearing were those whose completed $L$-functions
$\Lambda (s,f)$ have a minus sign in their functional equations.
Certainly this contains the subspace of forms whose central special
value vanishes, but there are additional forms that also contribute
(cf.~Example \ref{ex:127} below).
% Mark added `below' to avoid conflict with AGM2's Example 1.
 
Because of our extensive computations, we feel confident that
Conjecture \ref{conj:eisenstein} completely describes the Eisenstein
subspace for prime level.  However, Conjecture \ref{conj:eisenstein}
is not true for composite~$N$, as already remarked in the paragraph
after \cite[Example 1]{AGM2}.

\subsection{Paramodular forms}

\begin{conjecture}\label{conj:siegel}
For~$N$ prime, choose an equivalence class of eigenforms in
$P_{3}^{\ng } (N)$, and let $h$ be a representative.  Let $d_{h}$ be
the degree of the extension of $\Q$ generated by the eigenvalues of
$h$.  Then the cuspidal cohomology $H^{5}_{\cusp} (\Gamma_{0} (N);\C
)$ contains a $2d_{h}$-dimensional subspace spanned by Hecke
eigenclasses.  If $\xi$ is an eigenclass in this space, then up to
Galois conjugacy the Hecke polynomial $H (\xi )$ of $\xi$ from
\eqref{eqn:glheckpol} agrees with the Hecke polynomial $H_{\Sp} (h)$
of $h$ from \eqref{eqn:sympheckpol}.
\end{conjecture}

We remark that this equality means that $\xi $ is the functorial lift of $h$ with
respect to the natural inclusion of $L$-groups: ${}^L \GSp_4 \to
{}^L\GL_4$.

\subsection{Computational results} These are listed in Table
\ref{bettitab}, which shows our computed Betti numbers and the
dimensions of the constituents of the cohomology predicted by
Conjectures \ref{conj:eisenstein} and \ref{conj:siegel}.  For levels
$\leq 101$, we checked that the Hecke polynomial for $l = 2$ is
correct.

\begin{example}\label{ex:127}
We consider the case $N=127$.  There are two weight $2$ eigenforms,
with Hecke eigenvalues defining respectively a cubic and a septic
field.  There are three weight $4$ eigenforms, with Hecke eigenvalues
defining fields of degrees $1$, $13$, and $17$.  The degree $13$
eigenform has minus sign in the functional equation of its
$L$-function, which means its central special value vanishes.
However, there is also another vanishing at this level: the rational
eigenform also has vanishing central special value, vanishing that is
not forced by the sign of the functional equation.  Thus together
these modular forms account for a $2\times 10 + 14 = 34$ dimensional
subspace of $H^{5} (\Gamma_{0} (N); \C)$.
 
For the rest of the cohomology, we must consider $\SL_{3}$ and
paramodular contributions.  There is no cuspidal cohomology for
$\Gamma_{0}^{*} (127)\subset \SL_{3} (\Z)$.  The space of
non-Gritsenko lifts has dimension $3$.  Thus we see an additional
$6$-dimensional subspace of $H^{5}$ coming from these Siegel modular
forms, which means $\dim H^{5} (\Gamma_{0} (127); \C) \geq 40$.
Indeed our computations indicate that this Betti number equals $40$.
\end{example}

\begin{table}[ht]
\begin{center}
%\begin{tabular}{|c||c|c|c|c||c|}
\begin{tabular}{|p{50pt}||p{50pt}|p{50pt}|p{50pt}|p{50pt}||p{65pt}|}
\hline
Level&$S_{2} (N)$&$S_{4} (N)_{0}$&$\SL_{3}$&$P_{3}^{\ng} (N)$&$H^{5}
(\Gamma_{0} (N); \C )$\\
\hline
83&7&7&0&0&21\\
89&7&8&2&1&28\\
97&7&11&0&2&29\\
101&8&9&0&2&29\\
103&8&10&0&2&30\\
107&9&10&0&0&28\\
109&8&12&0&3&34\\
113&9&12&0&1&32\\
127&10&14&0&3&40\\
131&11&11&0&2&37\\
137&11&15&0&2&41\\
139&11&14&0&4&44\\
149&12&15&0&4&47\\
151&12&15&0&5&49\\
157&12&18&0&7&56\\
163&13&19&0&4&53\\
167&14&15&0&4&51\\
173&14&18&0&6&58\\
179&15&17&0&4&55\\
181&14&20&0&10&68\\
191&16&17&0&6&61\\
193&15&23&0&10&73\\
197&16&22&0&7&68\\
199&16&20&0&10&72\\
211&17&23&0&10&77\\
\hline
\end{tabular}
\end{center}
\medskip
\caption{Betti numbers and constituents of
Conjectures~\ref{conj:eisenstein}
and~\ref{conj:siegel}\label{bettitab}.  The entries are the dimensions
of the spaces in the headings, which are as follows: (i) $S_{2} (N)$
denotes weight $2$ cuspidal modular forms of level~$N$ and trivial
character, (ii) $S_{4} (N)_{0}$ denotes weight $4$ modular forms of
level~$N$, trivial character, and with vanishing central special
value, (iii) $\SL_{3}$ denotes the cuspidal cohomology of the
congruence subgroup $\Gamma_{0}^{*} (N)\subset \SL_{3} (\Z)$, and (iv)
$P^{\ng}_{3} (N)$ denotes weight $3$ paramodular forms of level~$N$
that are not Gritsenko lifts.  In all cases $2\times (\text{second} +
\text{fourth} + \text{fifth}) + \text{third}$ equals the dimension of
$H^{5} (\Gamma_{0} (N); \C)$.}
\end{table}

\begin{table}[ht]
\begin{center}
\begin{tabular}{|p{40pt}||p{80pt}|p{150pt}|}
\hline
\IIa & $\sigma_{2}\oplus \varepsilon^{2}\oplus \varepsilon^{3}$ & $(1-l^{2}T) (1-l^{3}T) (1-\alpha T + lT^{2})$\\
\hline
\IIb & $i\oplus \varepsilon^{2}\sigma_{2}\oplus \varepsilon  $ &$(1-T) (1-l T)(1-l^{2}\alpha T + l^{5}T^{2})$\\
\hline
\hline
\IV & $\sigma_{4}\oplus \varepsilon \oplus \varepsilon^{2}$ & $(1-lT) (1-l^{2}T) (1-\beta T + l^{3}T^{2})$\\
\hline
\hline
\IIIa & $\tau \oplus \varepsilon^{3} $&$(1-l^{3}T) (1-\gamma T + l \gamma 'T^{2} - l^{3}T^{3})$\\
\IIIb & $i\oplus \varepsilon \tau$&$(1-T) (1-l\gamma T + l^{3} \gamma 'T^{2} - l^{6}T^{3})$\\
\hline
\end{tabular}
\end{center}
\medskip
\caption{Galois representations and Hecke polynomials for Eisenstein
classes\label{grtab}.  See Sections \ref{s:eisenstein.paramodular} and~\ref{s:conjectures.computational} for explanation of notation.}
\end{table}

%%%%%%%%%%%%%%%%%%%%%%%%%%%%%%%%%%%%%%%%%%%%%%%%%%%%%%%%%%%%

\bibliographystyle{amsalpha_no_mr}
\bibliography{AGM-III}

\providecommand{\bysame}{\leavevmode\hbox to3em{\hrulefill}\thinspace}
\providecommand{\MR}{\relax\ifhmode\unskip\space\fi MR }
% \MRhref is called by the amsart/book/proc definition of \MR.
\providecommand{\MRhref}[2]{%
  \href{http://www.ams.org/mathscinet-getitem?mr=#1}{#2}
}
\providecommand{\href}[2]{#2}
\begin{thebibliography}{vGvdKTV97}

\bibitem[AGM02]{AGM1}
Avner Ash, Paul~E. Gunnells, and Mark McConnell, \emph{Cohomology of congruence
  subgroups of {${\rm SL}\sb 4(\Bbb Z)$}}, J. Number Theory \textbf{94} (2002),
  no.~1, 181--212.

\bibitem[AGM08]{AGM2}
\bysame, \emph{Cohomology of congruence subgroups of {${\rm SL}(4,\Bbb Z)$}.
  {II}}, J. Number Theory \textbf{128} (2008), no.~8, 2263--2274.

\bibitem[Ash92]{A}
Avner Ash, \emph{Galois representations attached to mod {$p$} cohomology of
  {${\rm GL}(n,{\bf Z})$}}, Duke Math. J. \textbf{65} (1992), no.~2, 235--255.

\bibitem[BCP97]{Magma}
Wieb Bosma, John Cannon, and Catherine Playoust, \emph{The {M}agma algebra
  system. {I}. {T}he user language}, J. Symbolic Comput. \textbf{24} (1997),
  no.~3-4, 235--265, Computational algebra and number theory (London, 1993).

\bibitem[BS73]{bs}
A.~Borel and J.-P. Serre, \emph{Corners and arithmetic groups}, Comm. Math.
  Helv. \textbf{48} (1973), 436--491.

\bibitem[Coh93]{Cohen}
Henri Cohen, \emph{A course in computational algebraic number theory}, Graduate
  Texts in Mathematics, vol. 138, Springer-Verlag, Berlin, 1993.

\bibitem[Dav06]{Davis}
Timothy~A. Davis, \emph{Direct methods for sparse linear systems}, Fundamentals
  of Algorithms, vol.~2, Society for Industrial and Applied Mathematics (SIAM),
  Philadelphia, PA, 2006.

\bibitem[DER89]{DER}
I.~S. Duff, A.~M. Erisman, and J.~K. Reid, \emph{Direct methods for sparse
  matrices}, second ed., Monographs on Numerical Analysis, The Clarendon Press
  Oxford University Press, New York, 1989, Oxford Science Publications.

\bibitem[DEVGU07]{DEvGU07}
Jean-Guillaume Dumas, Philippe Elbaz-Vincent, Pascal Giorgi, and Anna
  Urba{\'n}ska, \emph{Parallel computation of the rank of large sparse matrices
  from algebraic {$K$}-theory}, P{ASCO}'07, ACM, New York, 2007, pp.~43--52.

\bibitem[DSV01]{DSV}
Jean-Guillaume Dumas, B.~David Saunders, and Gilles Villard, \emph{On efficient
  sparse integer matrix {S}mith normal form computations}, J. Symbolic Comput.
  \textbf{32} (2001), no.~1-2, 71--99, Computer algebra and mechanized
  reasoning (St. Andrews, 2000).

\bibitem[EZ80]{ez}
M.~Eichler and D.~Zagier, \emph{Jacobi forms}, Birkh\"auser, 1980.

\bibitem[Gat00]{cll}
Erann Gat, \emph{{Lisp as an alternative to Java}}, Intelligence \textbf{11}
  (2000), no.~4, 21--24, \texttt{www.flownet.com/gat/papers/lisp-java.pdf}.

\bibitem[Gri95]{grit}
Valeri Gritsenko, \emph{Arithmetical lifting and its applications}, Number
  theory ({P}aris, 1992--1993), London Math. Soc. Lecture Note Ser., vol. 215,
  Cambridge Univ. Press, Cambridge, 1995, pp.~103--126.

\bibitem[GS]{M2}
Dan Grayson and Mike Stillman, \emph{Macaulay~2 software package},
  \texttt{www.math.uiuc.edu/Macaulay2}.

\bibitem[GVL96]{GVL}
Gene~H. Golub and Charles~F. Van~Loan, \emph{Matrix computations}, third ed.,
  Johns Hopkins Studies in the Mathematical Sciences, Johns Hopkins University
  Press, Baltimore, MD, 1996.

\bibitem[HHR93]{HHR93}
George Havas, Derek~F. Holt, and Sarah Rees, \emph{Recognizing badly presented
  {${\bf Z}$}-modules}, Linear Algebra Appl. \textbf{192} (1993), 137--163,
  Computational linear algebra in algebraic and related problems (Essen, 1992).

\bibitem[Ibu07]{ibuk2}
T.~Ibukiyama, \emph{{Dimension formulas of Siegel modular forms of weight $3$
  and supersingular abelian varieties}}, Proceedings of the Fourth Spring
  Conference on Modular Forms and Related Topics, Siegel Modular Forms and
  Abelian Varieties, 2007, pp.~39--60.

\bibitem[Jac85]{Jac}
Nathan Jacobson, \emph{{Basic Algebra I}}, 2nd ed., W. H. Freeman and Co.,
  1985.

\bibitem[LO91]{LaM-Odl}
B.~A. LaMacchia and A.~M. Odlyzko, \emph{Solving large sparse linear systems
  over finite fields}, Advances in Cryptology - CRYPTO '90, Lecture Notes in
  Computer Science, no. 537, Springer Verlag, 1991.

\bibitem[Mat03]{Matlab}
The MathWorks, \emph{{MATLAB}}, version 6.5.1, 2003.

\bibitem[McCa]{shh-lispwire}
Mark McConnell, \emph{\shh~2.1},
  \texttt{www.lispwire.com/entry-math-sheafhom-des}.

\bibitem[McCb]{shh}
\bysame, \emph{\shh\ software package},
  \texttt{www.geocities.com/mmcconnell17704/math.html}.

\bibitem[PS92]{PomS}
Carl Pomerance and J.~W. Smith, \emph{Reduction of huge, sparse matrices over
  finite fields via created catastrophes}, Experiment. Math. \textbf{1} (1992),
  no.~2, 89--94.

\bibitem[PY09]{py2}
C.~Poor and D.~Yuen, \emph{Paramodular cusp forms}, preprint, 2009.

\bibitem[Sch02]{sholten}
Jasper Scholten, \emph{{Mordell-Weil groups of elliptic surfaces and Galois
  representations}}, Ph.D. thesis, Groningen, 2002.

\bibitem[Tei98]{Tei}
J.~Teitelbaum, \emph{Euclid's algorithm and the {L}anczos method over finite
  fields}, Math. Comp. \textbf{67} (1998), no.~224, 1665--1678.

\bibitem[vdG08]{vdg}
Gerard van~der Geer, \emph{Siegel modular forms and their applications}, The
  1-2-3 of modular forms, Universitext, Springer, Berlin, 2008, pp.~181--245.

\bibitem[vGvdKTV97]{dutch}
Bert van Geemen, Wilberd van~der Kallen, Jaap Top, and Alain Verberkmoes,
  \emph{Hecke eigenforms in the cohomology of congruence subgroups of {${\rm
  SL}(3,\bold Z)$}}, Experiment. Math. \textbf{6} (1997), no.~2, 163--174.

\bibitem[Wie86]{Wie}
Douglas~H. Wiedemann, \emph{Solving sparse linear equations over finite
  fields}, IEEE Trans. Inform. Theory \textbf{32} (1986), no.~1, 54--62.

\end{thebibliography}

\end{document}